\documentclass[12pt]{article}
\usepackage{amscd}
\usepackage[all]{xy}
\usepackage{CJK}
\usepackage{CJKnumb}
\usepackage{graphicx}
\usepackage{amsfonts}
\usepackage{amsmath}
\usepackage{amssymb}
\usepackage{fancyhdr}
\usepackage{indentfirst}
\usepackage{titlesec}
\usepackage{mathrsfs}
\usepackage{bbm}
\usepackage{dsfont}
\usepackage{amsthm} %Ôö¼Ó¶¨ÀíÐλ·¾³£¬Ò²ÌṩÁËÖ¤Ã÷»·¾³
\usepackage{mathrsfs}

\newcommand{\cP}{\mathcal{P}}
\newcommand{\cM}{\mathcal{M}}

\newcommand{\cQ}{\mathcal{Q}}
\newcommand{\cR}{\mathcal{R}}

\newcommand{\cL}{\mathcal{L}}
\newcommand{\cN}{\mathcal{N}}
\newcommand{\cA}{\mathcal{A}}

\newcommand{\cC}{\mathcal{C}}

\newcommand{\cB}{\mathcal{B}}
\newcommand{\cK}{\mathcal{K}}
\newcommand{\cH}{\mathcal{H}}

\newcommand{\cS}{\mathcal{S}}

\theoremstyle{definition}  \newtheorem{Def}{Definition}
\theoremstyle{plain}  \newtheorem{Thm}{Theorem} \theoremstyle{plain}
\newtheorem{cor}{Corollary}
 \theoremstyle{remark} \newtheorem{Rek}{Remark}
 \theoremstyle{plain}\newtheorem{Lem}{Lemma}
 \theoremstyle{plain}
\newtheorem{Prop}{Proposition}
\theoremstyle{remark}

\begin{document}

\title{\textbf{Higher Dimensional Homology Algebra IV:Projective Resolutions and Derived 2-Functors in ($\cR$-2-Mod)}}
\author{Fang HUANG, Shao-Han CHEN, Wei CHEN, Zhu-Jun ZHENG\thanks{Supported in part by NSFC with grant Number
10971071
 and Provincial Foundation of Innovative
Scholars of Henan.} }
\date{}
 \maketitle

\begin{center}
\begin{minipage}{5in}
{\bf  Abstract}: In this paper, we will construct the projective
resolution of any $\cR$-2-module, define the derived 2-functor and
give some related properties of the derived 2-functor.

\\
%\bf CLC number: O154.1
\end{minipage}
\end{center}
\maketitle \hspace{1cm}
%\begin{center}
%\Large\bf Á½ÖÖ 2-Ä£µÄ¹Øϵ\bigbreak \normalsize\it{\Large\bf
%»Æ·¼$^{1,}$\footnote{{\bf Biography}: Fang Huang(1984-), female,
%%native of Shangqiu, Henan, Ph.D., engages in algebraic geometry,
%category theory. E-mail:
%huang.fang@mail.scut.edu.cn}\bigbreak
% (1. »ªÄÏÀí¹¤´óѧÀíѧԺÊýѧϵ£¬¹ãÖÝ 510640 £»
%2. ºÓÄÏ´óѧ£¬¿ª·â 475004)
%\end{center}
%\hspace{1cm}
%\begin{center}
%\begin{minipage}{5in}
%{\bf  Õª\quad Òª}: Ä¿Ç°Ö÷ÒªÓÐÁ½Àà¹ØÓÚÄ£µÄ·¶³ë»¯¸ÐÄһÀàÊÇ
%\cite{1}ÖеÄ2-Ä££¬ÁíÒ»ÀàÊÇ\cite{6} ÖеÄ$\cR$-2-Ä£.
%±¾ÎÄÖ÷ÒªÖ¤Ã÷ÕâÁ½ÖÖ¶¨ÒåµÄµÈ¼ÛÐÔ£¬
%ÕâÖֵȼÛÐÔÔÚÑо¿¸ßά´úÊýÀíÂÛÖÐÓÐÖØÒªÓá£
%\\
%{\bf{¹Ø¼ü´Ê}:} ·¶³ë»¯; 2-Ä£; $\cR$2-Ä£\\
%\end{minipage}
%\end{center}
%\maketitle

\section{Introduction}

A 2-ring $\cR$ is a category with categorical ringed structure(see
\cite{5}). As 1-dimensional algebra, we defined $\cR$-2-modules
\cite{4} in a different way with M.Dupont's 2-modules in his PhD.
thesis\cite {2}. An $\cR$-2-module we mentioned in this paper is
$(\cA,I,\cdot,a,b,i,z)$, where $\cA$ is a symmetric 2-group with
$\cR$-2-module structure $\cdot$, $I$ is the unit object under
$\cdot$, $a,b,i,z$ are natural isomorphisms satisfying canonical
properties \cite{4}.

%In \cite{11}, A.del R\'{\i}o, J. Mart\'{\i}nez-Moreno  and E. M.
%Vitale gave the definition of cohomology symmetric 2-groups for any
%complex of symmetric 2-groups in the 2-category (2-SGp)(which is an
%abelian 2-category\cite{2}) of symmetric categorical groups(we call
%them symmetric 2-groups) after discussing the relative kernel and
%the relative cokernel, and constructed a long 2-exact sequence from
%an extension of complexes. %
Based on the works of A.del R\'{\i}o, J. Mart\'{\i}nez-Moreno and E.
M. Vitale\cite{11}, we defined the left derived 2-functor in the
2-category (2-SGp) and gave a fundamental property of the derived
2-functor in our third paper \cite{20} of the series of higher
dimensional homology algebra. In \cite{2,4}, the authors showed that
the 2-category ($\cR$-2-Mod) is an abelian 2-category which has
enough projective(injective) objects(\cite{14,21}). Naturally, we
will consider the higher dimensional homological theory in
($\cR$-2-Mod).

The aim of this paper is to develop a homological theory in the
2-category ($\cR$-2-Mod) just like the 1-dimensional case. We will
construct the projective resolution of any $\cR$-2-module, which is
unique up to 2-chain homotopy(Definition 3) and give the definition
of the left derived 2-functor in ($\cR$-2-Mod). Moreover, we shall
give a fundamental property of the derived 2-functor. In our paper,
most results are similar to \cite{20}, just replacing the morphisms
of symmetric 2-groups by morphisms of $\cR$-2-modules. The most
different and difficult
%between this paper and
%the last paper(\cite{20})
are to give the $\cR$-2-module structures of relative kernel and
cokernel.

The present paper is organized as follows. In section 2, we give
some basic facts on $\cR$-2-modules such as the relative kernel and
cokernel which are appeared in \cite{2,11,6} for symmetric 2-group
case. The homology $\cR$-2-modules of a complex of $\cR$-2-modules
appear in this section, too. In the last section, we mainly give the
definition of projective resolution of an $\cR$-2-module and give
its construction(Proposition 2). After the basic definition of
derived 2-functors from abelian 2-category ($\cR$-2-Mod) to
($\cS$-2-Mod)(\cite{2,4}), we obtain our main result(Theorem 2).

This is the fourth paper of the series works on higher dimensional
homological algebra.
%ollowed by our last paper(\cite{20})
%The first paper is "2-Modules and the Representation of
%2-Rings\cite{4}". The second paper is "Higher Dimensional Homology
%Algebra II: Projectivity\cite{14}". The recent paper is "Higher
%Dimensional Homology Algebra III:Projective Resolutions and Derived
%2-Functors in ($\cR$-2-Mod)". These papers form a basic theory of
%higher dimensional homological algebra, we will try to complete this
%whole theory in our researching works.

\section{Preliminary}
In this section, we give the definitions and constructions of the
relative (co)kernel
%and relative 2-exactness of a sequence
 in ($\cR$-2-Mod) from the
definitions of them given in \cite{2,11}, and then give the homology
$\cR$-2-modules of a complex of $\cR$-2-modules which is similar to
the homology symmetric 2-groups given in \cite{20}, where $\cR$ is a
2-ring. In this paper, we will omit the composition symbol $\circ$
in our diagrams.
%\begin{Def}\cite{1,11,6} For a sequence
%$(F,\varphi,G)$ in (2-SGp) as in the following diagram:
%\begin{center}
%\scalebox{0.9}[0.85]{\includegraphics{p1.eps}}
%{\footnotesize Fig.1.}
%\end{center}
%By the universal properties of kernel and cokernel(\cite{2,11,6}),
%there are homomorphisms $F_0,G_0$ as in the following diagram:
%\begin{center}
%\scalebox{0.9}[0.85]{\includegraphics{p2.eps}}
%{\footnotesize Fig.2.}
%\end{center}
%The sequence $(F,\varphi,G)$ is 2-exact if it satisfies one of the
%following equivalent conditions:

%1) $F_0:\cA\rightarrow KerG$ is full and essentially surjective;

%2) $G_0:CokerF\rightarrow \cC$ is full and faithful.
%\end{Def}

%\begin{Rek} There are four equivalent conditions in above
%definition from Proposition 6.2 in \cite{6}.
%\end{Rek}
\begin{Def}
The relative kernel of the sequence
$(F,\varphi,G):\cA\rightarrow\cB\rightarrow\cC$ in ($\cR$-2-Mod)is
the triple
$(Ker(F,\varphi),e_{(F,\varphi)},\varepsilon_{(F,\varphi)})$ in
($\cR$-2-Mod) as in the following diagram
\begin{center}
\scalebox{0.9}[0.85]{\includegraphics{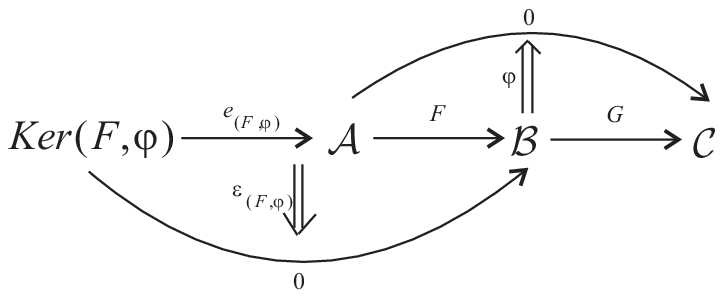}}
%{\footnotesize Fig.1.}
\end{center}
with $\varepsilon_{(F,\varphi)}$ compatible with $\varphi$, i.e. the
following diagram commutes
\begin{center}
\scalebox{0.9}[0.85]{\includegraphics{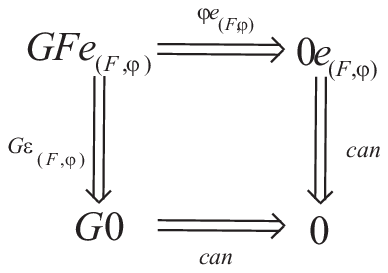}}
%{\footnotesize Fig.1.}
\end{center}
and satisfies the following universal property:

Given a diagram in ($\cR$-2-Mod)
\begin{center}
\scalebox{0.9}[0.85]{\includegraphics{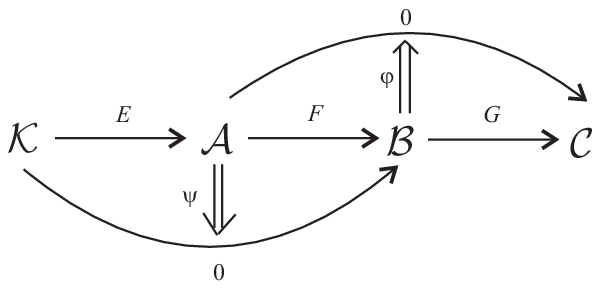}}
%{\footnotesize Fig.1.}
\end{center}
with $\psi$ compatible with $\varphi$, there is a factorization
$$
(E^{'}:\cK\rightarrow Ker(F,\varphi),\psi^{'}:e_{(F,\varphi)}\circ
E^{'}\Rightarrow E)
$$
in ($\cR$-2-Mod) through
$(e_{(F,\varphi)},\varepsilon_{(F,\varphi)})$, that is the following
diagram commutes
\begin{center}
\scalebox{0.9}[0.85]{\includegraphics{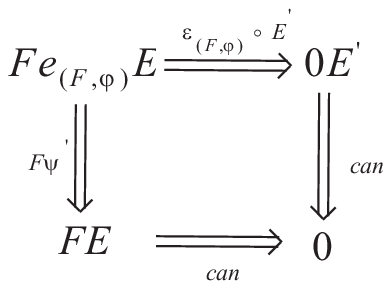}}
%{\footnotesize Fig.1.}
\end{center}
and if $(E^{''},\psi^{''})$ is another factorization of $(E,\psi)$
through $(e_{(F,\varphi)},\varepsilon_{(F,\varphi)})$, then there is
a unique 2-morphism $e:E^{'}\Rightarrow E^{''}$, such that
\begin{center}
\scalebox{0.9}[0.85]{\includegraphics{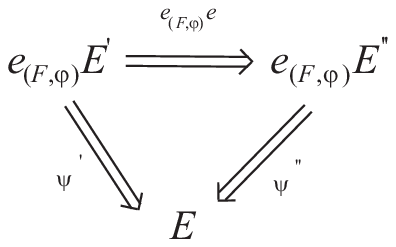}}
%{\footnotesize Fig.1.}
\end{center}
commutes.
\end{Def}

The existence of relative kernel is given similarly to the general
kernel\cite{4}.

First, $Ker(F,\varphi)$ is a symmetric 2-group(see \cite{11,20})
consisting of:

$\cdot$ An object is a pair $(A\in obj(\cA), a:F(A)\rightarrow 0)$
such that the following diagram commutes
\begin{center}
\scalebox{0.9}[0.85]{\includegraphics{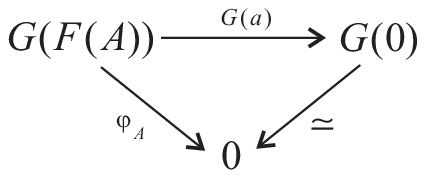}}
\end{center}

$\cdot$ A morphism $f:(A,a)\rightarrow(A^{'},a^{'})$ is a morphism $
f:A\rightarrow A^{'}$ of $\cA$ such that the following diagram
commutes
\begin{center}
\scalebox{0.9}[0.85]{\includegraphics{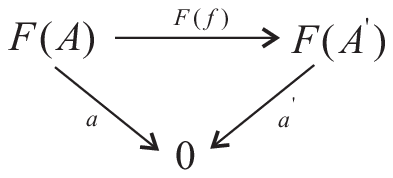}}
\end{center}

Second, $Ker(F,\varphi)$ is an $\cR$-2-module:

There is a bifunctor
\begin{align*}
&\cdot: \cR\times Ker(F,\varphi)\rightarrow Ker(F,\varphi)\\
&\hspace{1.5cm}(r,(A,a))\mapsto r\cdot(A,a)\triangleq (r\cdot A, r\cdot a),\\
&\hspace{2.3cm}(x,f)\mapsto x\cdot f,
%(r_1,(A_1,a_1))\xrightarrow[]{(x,f)}(r_2,(A_2,a_2))\mapsto
%r\cdot(A)
\end{align*}
where $r\cdot A$ and $x\cdot f$ are the object and morphism in $\cA$
under its $\cR$-2-module structure, respectively, $r\cdot a$ is the
composition morphism $F(r\cdot A)\backsimeq r\cdot
F(A)\xrightarrow[]{r\cdot a}r\cdot 0\backsimeq 0$. The above
bifunctor is well-defined. In fact, for $(A,a)\in
obj(Ker(F,\varphi))$ with $G(a)=\varphi_{A}$, there is $G(r\cdot
a)=r\cdot G(a)=r\cdot \varphi_{A}=\varphi_{r\cdot A}$ from the basic
properties of $\cR$-2-modules. Moreover, the natural isomorphisms in
the definition of $\cR$-2-modules and the universal property are
given as general kernels(more details see \cite{4}).

%$\cdot$ The faithful functor
%$e_{(F,\varphi)}:Ker(F,\varphi)\rightarrow\cA$ is defined by
%$e_{(F,\varphi)}(A,a)=A$, and the natural transformation
%$\varepsilon_{(F,\varphi)}:F\circ e_{(F,\varphi)}\Rightarrow 0$ by
%$(\varepsilon_{(F,\varphi)})_{(A,a)}=a$.

\begin{Def}
The relative cokernel of the sequence
$(F,\varphi,G):\cA\rightarrow\cB\rightarrow\cC$ in ($\cR$-2-Mod) is
the triple $(Coker(\varphi,G),p_{(\varphi,G)},\pi_{(\varphi,G)})$ in
($\cR$-2-Mod) as in the following diagram
\begin{center}
\scalebox{0.9}[0.85]{\includegraphics{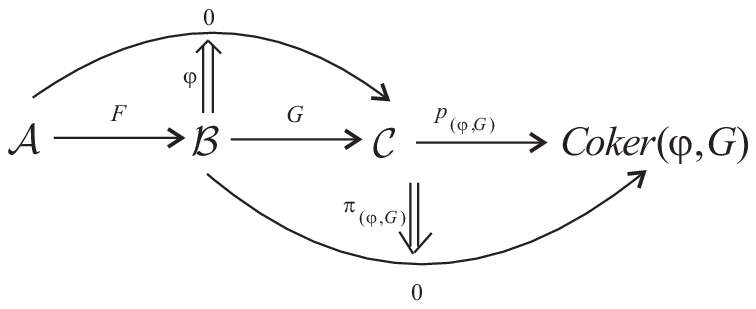}}
%{\footnotesize Fig.1.}
\end{center}
with $\pi_{(\varphi,G)}$ compatible with $\varphi$, i.e. the
following diagram commutes
\begin{center}
\scalebox{0.9}[0.85]{\includegraphics{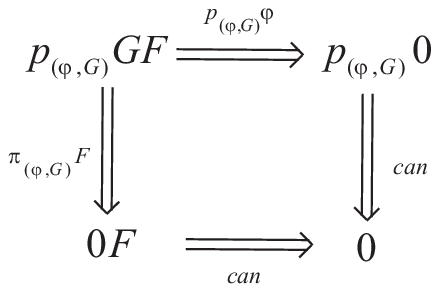}}
%{\footnotesize Fig.1.}
\end{center}
and satisfies the following universal property:

Given a diagram in ($\cR$-2-Mod)
\begin{center}
\scalebox{0.9}[0.85]{\includegraphics{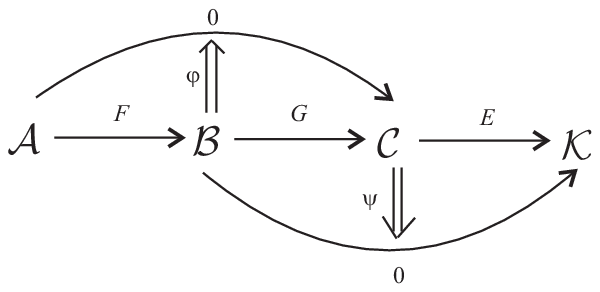}}
%{\footnotesize Fig.1.}
\end{center}
with $\psi$ compatible with $\varphi$, there is a factorization
$$
(E^{'}: Coker(\varphi,G)\rightarrow\cK,\psi^{'}: E^{'}\circ
p_{(\varphi,G)}\Rightarrow E)
$$
in ($\cR$-2-Mod) through $(p_{(\varphi,G)},\pi_{(\varphi,G)})$, that
is the following diagram commutes
\begin{center}
\scalebox{0.9}[0.85]{\includegraphics{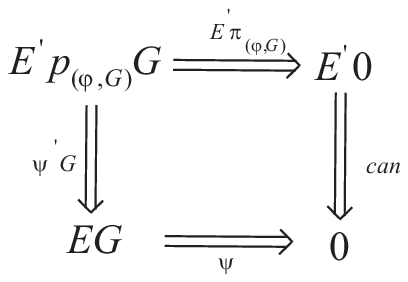}}
%{\footnotesize Fig.1.}
\end{center}
and if $(E^{''},\psi^{''})$ is another factorization of $(E,\psi)$
through $(p_{(\varphi,G)},\pi_{(\varphi,G})$, then there is a unique
2-morphism $e:E^{'}\Rightarrow E^{''}$, such that
\begin{center}
\scalebox{0.9}[0.85]{\includegraphics{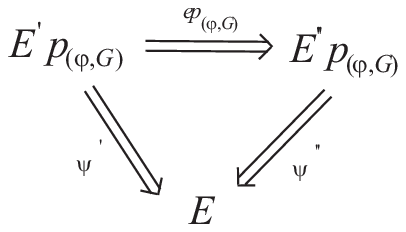}}
%{\footnotesize Fig.1.}
\end{center}
commutes.
\end{Def}

The existence of relative cokernel is also similar to the general
cokernel\cite{4}.

First, $Coker(\varphi,G)$ is a symmetric 2-group(see \cite{11,20})
consisting of:

$\cdot$  Objects are those of $\cC$.

$\cdot$ A morphism from $X$ to $Y$ is an equivalent class of a pair
$(B,f):X\rightarrow Y$ with $B\in obj(\cB)$ and $f:X\rightarrow
G(B)+Y$. For two morphisms $(B,f),(B^{'},f^{'}):X\rightarrow Y$ are
equivalent if there is $A\in obj(\cA)$ and $a:B\rightarrow
F(A)+B^{'}$ such that the following diagram commutes
\begin{center}
\scalebox{0.9}[0.85]{\includegraphics{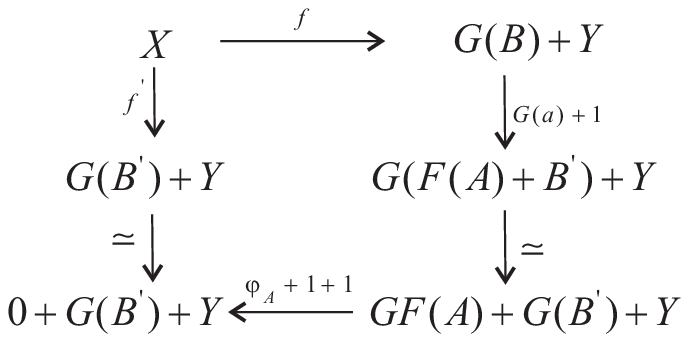}}
%{\footnotesize Fig.1.}
\end{center}

Second, $Coker(\varphi,G)$ is an $\cR$-2-module.

There is a bifunctor
\begin{align*}
&\hspace{0.5cm}\cdot: \cR\times Coker(\varphi,G)\rightarrow Coker(\varphi,G)\\
&\hspace{3.1cm}(r,X)\mapsto r\cdot X,\\
&(r_1\xrightarrow[]{x}r_2,X\xrightarrow[]{[B,f]}Y)\mapsto r_{1}\cdot
X\xrightarrow[]{[r_{1}\cdot B,\overline{f}]}r_{2}\cdot Y,
%(r_1,(A_1,a_1))\xrightarrow[]{(x,f)}(r_2,(A_2,a_2))\mapsto
%r\cdot(A)
\end{align*}
where $r\cdot X$ and $r_1\cdot B$ are the objects in $\cC$ and $\cB$
under the $\cR$-2-module structures of them,respectively,
$\overline{f}$ is the composition morphism $r_1\cdot
X\xrightarrow[]{r_1\cdot f} r_1\cdot(G(B)+Y)\backsimeq r_1\cdot
G(B)+r_1\cdot Y\backsimeq G(r_1\cdot B)+r_1\cdot
Y\xrightarrow[]{1+x\cdot Y}G(r_1\cdot B)+r_2\cdot Y$. This bifunctor
is well-defined. In fact, if $[B,f]=[B^{'},f^{'}]:X\rightarrow Y$,
i.e. there exist $A\in obj(\cA)$ and $a:B\rightarrow F(A)+B^{'}$
such that the following diagram commutes
\begin{center}
\scalebox{0.9}[0.85]{\includegraphics{p9.3.eps}}
%{\footnotesize Fig.1.}
\end{center}
Hence, there exist $r_1\cdot A\in obj(\cA)$ and $r_1\cdot a:r_1\cdot
B\rightarrow F(r_1\cdot A)+r_1\cdot B^{'}$ such that the following
diagram commutes
\begin{center}
\scalebox{0.8}[0.8]{\includegraphics{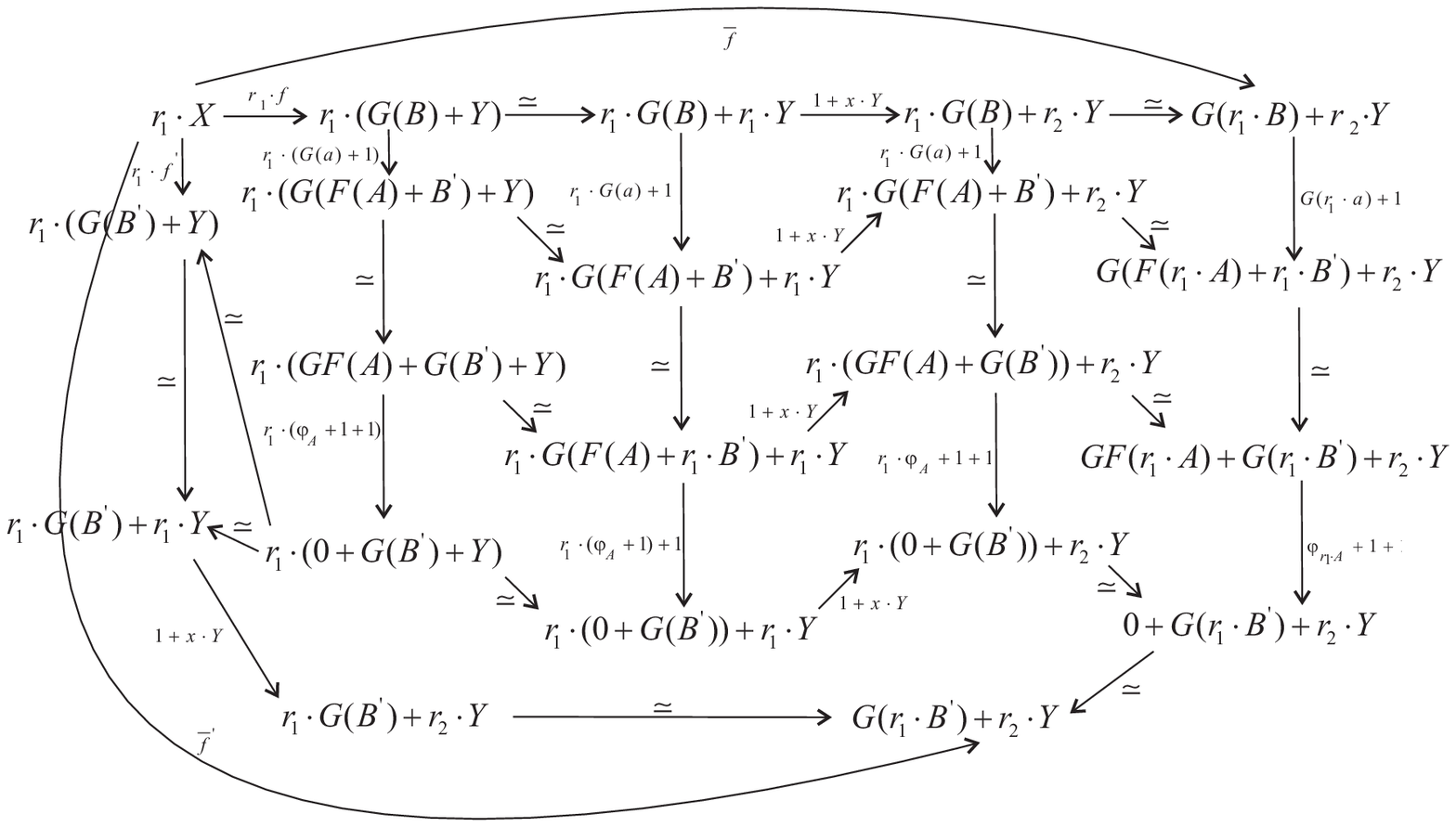}}
%{\footnotesize Fig.1.}
\end{center}
Then $x\cdot[B,f]=[r_1\cdot B,\overline{f}]=[r_1\cdot
B^{'},\overline{f}^{'}]=r_1\cdot[B^{'},f^{'}]:r_1\cdot X\rightarrow
r_2\cdot Y$.

Moreover, the natural isomorphisms in the definitions of
$\cR$-2-module and the universal property are given as general
cokernel(more details see \cite{4}).

%$\cdot$ The essentially surjective functor
%$p_{(\varphi,G)}:\cC\rightarrow Coker(\varphi,G)$ is defined by
%$p_{(\varphi,G)}(X)=X$, and the natural transformation
%$\pi_{(\varphi,G)}: p_{(\varphi,G)}\circ G\Rightarrow 0$ by
%$(\pi_{(\varphi,G)})_{B}=1_{G(B)}$.

%The universal properties of relative kernel and cokernel just like
%the usual ones, more details see \cite{11}.

%\begin{Def}
%Consider the following diagram in ($\cR$-2-Mod)
%\begin{center}
%\scalebox{0.9}[0.85]{\includegraphics{p10.eps}}
%\end{center}
%with $\alpha$ compatible with $\varphi$ and $\varphi$ compatible
%with $\gamma$. By the universal property of the relative kernel
%$Ker(G,\gamma)$, we get a factorization $(F^{'},\varphi^{'})$ of
%$(F,\varphi)$ through $(e_{(F,\varphi)},\varepsilon_{(F,\varphi)})$.
%By the cancelation property of $e_{(F,\varphi)}$, we have a
%2-morphism $\overline{\alpha}$ as in the following diagram
%\begin{center}
%\scalebox{0.9}[0.85]{\includegraphics{p11.eps}}
%{\footnotesize Fig.1.}
%\end{center}
%We say that the sequence $(L,\alpha,F,\varphi,G,\gamma,M)$ is
%relative 2-exact in $\cB$ if the $\cR$-homomorphism $F^{'}$ is
%essentially surjective and $\overline{\alpha}$-full(\cite{11}).
%\end{Def}
%\begin{Rek} The equivalent definition of relative 2-exact can also
%be given similar as this definition for symmetric 2-groups in
%\cite{11}.
%\end{Rek}
\begin{Rek}Just like symmetric 2-group discussed in \cite{1,2,11,13,20},
we can give the definitions of (relative-)2-exact, cohomology
$\cR$-2-modules in ($\cR$-2-Mod).
\end{Rek}

A complex of $\cR$-2-modules in ($\cR$-2-Mod) is a sequence
$$
\cA_{\cdot}=\cdot\cdot\cdot\xrightarrow[]{L_{n+1}}\cA_{n}\xrightarrow[]
{L_{n}}\cA_{n-1}\xrightarrow[]{L_{n-1}}\cA_{n-2}\xrightarrow[]{L_{n-2}}\cdot\cdot\cdot
\xrightarrow[]{L_{2}}\cA_{1}\xrightarrow[]{L_{1}}\cA_{0}
$$
together with a family of 2-morphisms $\{\alpha_{n}:L_{n-1}\circ
L_{n}\Rightarrow 0\}_{n\geq2}$ in ($\cR$-2-Mod) such that, for all
$n$, the following diagram commutes
\begin{center}
\scalebox{0.9}[0.85]{\includegraphics{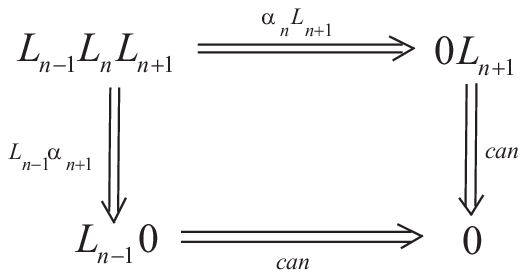}}
\end{center}

Consider part of the complex
\begin{center}
\scalebox{0.9}[0.85]{\includegraphics{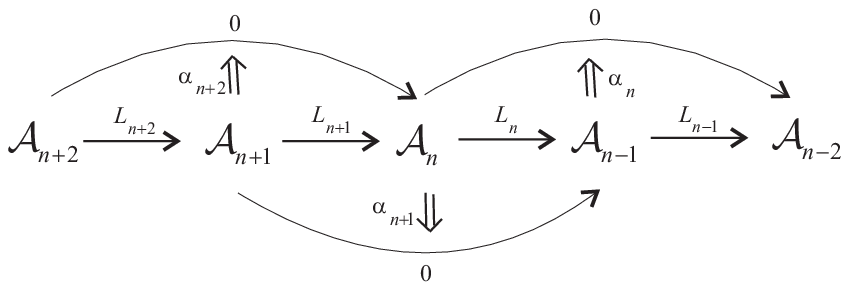}}
\end{center}
Based on the universal properties of relative kernel
$Ker(L_n,\alpha_n)$, we have the following diagram
\begin{center}
\scalebox{0.9}[0.85]{\includegraphics{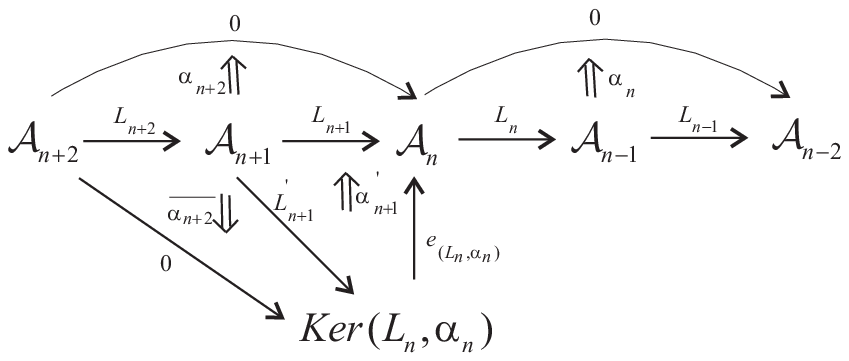}}
\end{center}
Similarly as the definition of (co)homology 2-group in \cite{11,20},
the $n$th homology $\cR$-2-module $\cH_{n}(\cA_{\cdot})$ of the
complex $\cA_{\cdot}$ defined as the relative cokernel
$Coker(\overline{\alpha_{n+2}},L_{n+1}^{'})$.

Note that, to get $\cH_{0}(\cA_{\cdot})$ and $\cH_{1}(\cA_{\cdot})$,
we have to complete the complex $\cA_{\cdot}$ on the right with the
two zero morphisms and two canonical 2-morphisms\\
$ \cdot\cdot\cdot\xrightarrow
{L_2}\cA_1\xrightarrow[]{L_1}\cA_{0}\xrightarrow[]{0}0\xrightarrow[]{0}0,
$ can : $0\circ L_1\Rightarrow 0,$ can : $0\circ 0\Rightarrow 0.$

The explicit description of $\cH_{n}(\cA_{\cdot})$ can also be given
from the existence of relative kernel and relative cokernel in
($\cR$-2-Mod) like the symmetric 2-group case in \cite{11,20}.

%We give an explicit description of $\cH_{n}(\cA_{\cdot})$ following
%from the cohomology symmetric 2-group given in \cite{11}.
%
%$\cdot$ an object of $\cH_{n}(\cA_{\cdot})$ is an object of the
%relative kernel $Ker(L_{n},\alpha_n)$, that is a pair
%$$
%(A_n\in obj(\cA_n), a_n:L_n(A_n)\rightarrow 0)
%$$
%such that $L_{n-1}(a_n)=(\alpha_{n})_{A_n}$;
%
%$\cdot$ a morphism $(A_n,a_n)\rightarrow (A_n^{'},a_n^{'})$ is an
%equivalent pair
%$$
%(X_{n+1}\in obj(\cA_{n+1}),x_{n+1}:A_n\rightarrow
%L_{n+1}(X_{n+1})+A_n^{'})
%$$
%such that the following diagram commutes
%\begin{center}
%\scalebox{0.9}[0.85]{\includegraphics{p11.1.eps}}
%\end{center}
%For two morphisms
%$(X_{n+1},x_{n+1}),(X_{n+1}^{'},x_{n+1}^{'}):(A_n,a_n)\rightarrow
%(A_n^{'},a_n^{'})$ are equivalent if there is a pair
%$$
%(X_{n+2}\in obj(\cA_{n+2}),x_{n+2}:X_{n+1}\rightarrow
%L_{n+2}(X_{n+2})+A_{n+1}^{'})
%$$
%such that the following diagram commutes
%\begin{center}
%^\scalebox{0.9}[0.85]{\includegraphics{p11.2.eps}}
%\end{center}

%Similarly \cite{11},
A morphism
$(F_{\cdot},\lambda_{\cdot}):\cA_{\cdot}\rightarrow\cB_{\cdot}$ of
complexes in ($\cR$-2-Mod) is a picture in the following diagram
such that the following diagram commutes
\begin{center}
\scalebox{0.9}[0.85]{\includegraphics{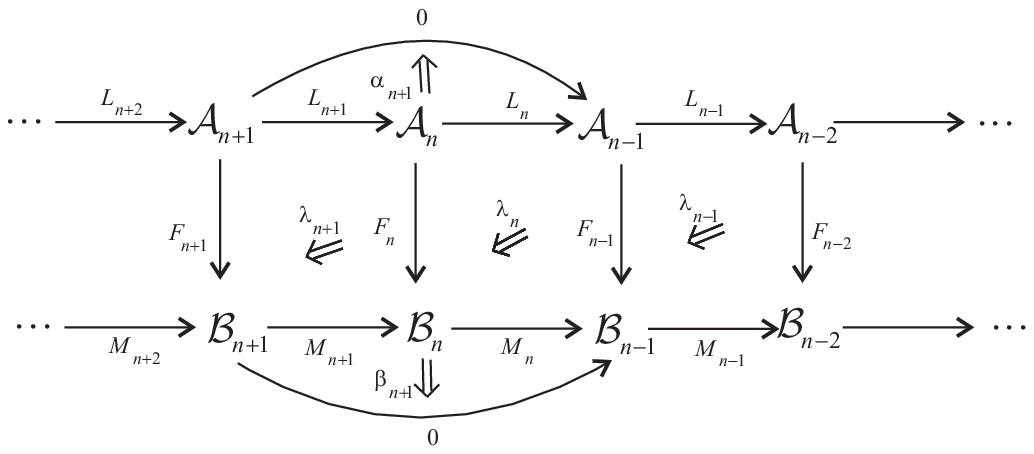}}
\end{center}
where $F_{n}:\cA_n\rightarrow \cB_n$ is 1-morphism in ($\cR$-2-Mod),
$\lambda_{n}:F_{n-1}\circ L_{n}\Rightarrow M_{n}\circ F_n$ is
2-morphism in ($\cR$-2-Mod), for each $n$, making the following
diagram commutative
\begin{center}
\scalebox{0.9}[0.85]{\includegraphics{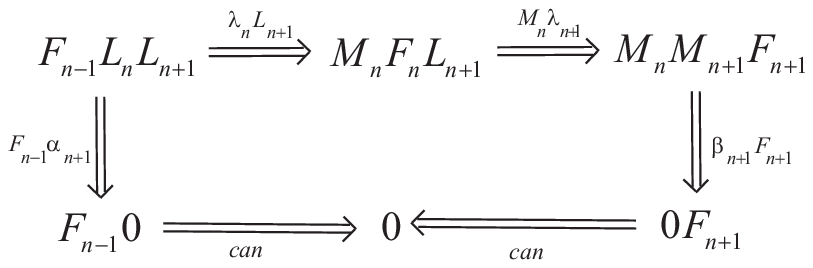}}
% {\footnotesize
%Fig.1.}
\end{center}
Such a morphism induces, for each $n$, a morphism of homology
$\cR$-2-modules $\cH_{n}(F_{\cdot}):\cH_{n}(\cA_{\cdot})\rightarrow
\cH_{n}(\cB_{\cdot})$ from the universal properties of relative
kernels and cokernels. It can be described as follows:(more details
see \cite{11,20} ).

Given an object $(A_n\in obj(\cA_{n}),a_n:L_n(A_n)\rightarrow 0)$ of
$\cH_{n}(\cA_{\cdot})$,
%with $L_{n-1}(a_n)=(\alpha_n)_{A_n}$,
we have $\cH_{n}(F_{\cdot})(A_n,a_n)=(F_n(A_n)\in
obj(\cB_{n}),b_n:M_{n}(F_n(A_n))\rightarrow 0)$, where $b_n$ is the
composition
$M_{n}(F_n(A_n))\xrightarrow[]{(\lambda_n)_{A_n}^{-1}}F_{n-1}L_{n}(A_n)\xrightarrow[]{F_{n-1}(a_n)}F_{n-1}(0)\backsimeq
0.$
%, together with $M_{n-1}(b_n)=(\beta_n)_{F_n(A_n)}$.
%In fact, from
%the commutative diagram of $\lambda_n$, we have the following
%commutative diagram
%\begin{center}
%\scalebox{0.9}[0.85]{\includegraphics{p11.5.eps}}
% {\footnotesize
%Fig.1.}
%\end{center}
%Moreover, consider 2-morphism $\lambda_{n-1}:F_{n-2}\circ
%L_{n-1}\Rightarrow M_{n-1}\circ F_{n-1}$ and a morphism
%$a_n:L_n(A_n)\rightarrow 0$ in $\cA_{n-1}$, we have the following
%commutative diagram
%\begin{center} \scalebox{0.9}[0.85]{\includegraphics{p11.6.eps}}
% {\footnotesize
%Fig.1.}
%\end{center}
%Then, by the above two diagrams, we have
%$M_{n-1}(b_n)=(\beta_n)_{F_n(A_n)}$, i.e. $(F_n(A_n),b_n)$ is an
%object of $\cH_{n}(\cB_{\cdot})$.

Given a morphism $[X_{n+1}\in obj(\cA_{n+1}),x_{n+1}:A_n\rightarrow
L_{n-1}(X_{n+1})+A_{n}^{'}]:(A_n,a_n)\rightarrow(A_{n}^{'},a_{n}^{'})$
in $\cH_{n}(\cA_{\cdot}),$
%satisfying the condition as in the above definition.
we have $ \cH_{n}(F_{\cdot})[X_{n+1},x_{n+1}]=[F_{n+1}(X_{n+1})\in
obj(\cB_{n+1}),\overline{x_{n+1}}:F_{n}(A_n)\rightarrow
M_{n+1}(F_{n+1}(X_{n+1}))+F_{n}(A_{n}^{'})]:(F_{n}(A_{n},b_n)\rightarrow(F_n(A_{n}^{'}),b_{n}^{'}),
$ where $\overline{x_{n+1}}$ is the composition
$F_{n}(A_n)\xrightarrow[]{F_{n}(x_{n+1})}F_{n}(L_{n+1}(X_{n+1}+A_{n}^{'}))\backsimeq
F_{n}L_{n+1}(X_{n+1})+F_n(A_{n}^{'})\xrightarrow[]{(\lambda_{n+1})_{X_{n+1}}+1}M_{n+1}(F_{n+1}(X_{n+1}))+F_{n}(A_{n}^{'})
$.
% and such that the following diagram commutes
%\begin{center}
%\scalebox{0.9}[0.85]{\includegraphics{p11.7.eps}}
% {\footnotesize
%Fig.1.}
%\end{center}
%In fact, we have the following commutative diagrams
%\begin{center} \scalebox{0.9}[0.85]{\includegraphics{p11.8.eps}}
% {\footnotesize
%Fig.1.}
%\end{center}
%The commutativity of I follows from $\lambda_{n}$ is a natural
%transformation. The commutativity of II follows from $\lambda_{n}$
%is a 2-morphism. The commutativity of III follows from the
%commutativity of $\lambda_{n}$ in definition. The commutativity of
%IV is obvious. The commutativity of V follows from the operation of
%$F_{n-1}$ on the commutative diagram of $[X_{n+1},x_{n+1}]$.

%$\cH_{n}(F_{\cdot})$ is a morphism in ($\cR$-2-Mod) follows from the
%properties of $F_{n}$.

%For a morphism $F_{\cdot}:\cA_{\cdot}\rightarrow \cB_{\cdot}$ of
%complexes in ($\cR$-2-Mod). There is an induced morphism
%$\cH_{n}(F_{\cdot}):\cH_{n}(\cA_{\cdot})\rightarrow
%\cH_{n}(\cB_{\cdot})$ of homology $\cR$-2-modules ( \cite{11,20}
%about symmetric 2-groups).

\begin{Rek}
1. For a complex of $\cR$-2-modules which is relative 2-exact in
each point, the (co)homology $\cR$-2-modules are always zero
$\cR$-2-module(only one object and one morphism).

2. For morphisms
$\cA_{\cdot}\xrightarrow[]{(F_{\cdot},\lambda_{\cdot})}\cB_{\cdot}\xrightarrow[]{(G_{\cdot},\mu_{\cdot})}\cC_{\cdot}$
of complexes in ($\cR$-2-Mod), their composition is given by
$(G_{n}\circ F_{n},(\mu_n\circ F_{n+1})\star (G_n\circ \lambda_n))$,
for $n\in \mathds{Z}$, where $\star$ is the vertical composition of
2-morphisms in 2-category(\cite{19}). Moreover,
$\cH_{n}(G_{\cdot}\circ F_{\cdot})\backsimeq \cH_{n}(G_{\cdot})\circ
\cH_n(F_{\cdot}).$
%of homology $\cR$-2-modules.
\end{Rek}

\begin{Def}
Let $(F_{\cdot},\lambda_{\cdot}),
(G_{\cdot},\mu_{\cdot}):(\cA_{\cdot},L_{\cdot},\alpha_{\cdot})\rightarrow
(\cB_{\cdot},M_{\cdot},\beta_{\cdot})$ be two morphisms of complexes
of $\cR$-2-modules. If there is a family of 1-morphisms
$\{H_{n}:\cA_{n}\rightarrow\cB_{n+1}\}_{n\in \mathds{Z}}$ and a
family of 2-morphisms $\{\tau_{n}:F_{n}\Rightarrow M_{n+1}\circ
H_{n}+H_{n-1}\circ L_n+G_{n}:\cA_{n}\rightarrow \cB_n\}_{n\in
\mathds{Z}}$ satisfying the obvious compatible conditions, i.e. the
following diagram commutes
\begin{center}
\scalebox{0.9}[0.85]{\includegraphics{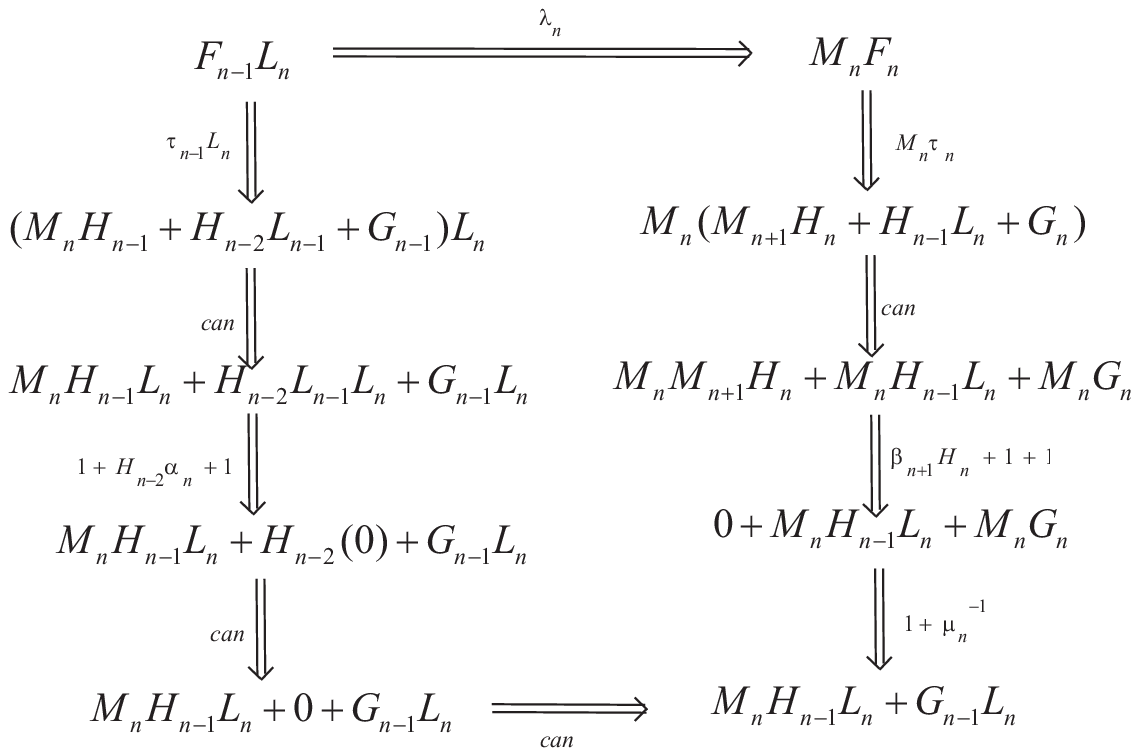}}
\end{center}
We call the above morphisms $(F_{\cdot},\lambda_{\cdot}),
(G_{\cdot},\mu_{\cdot})$ are 2-chain homotopy in ($\cR$-2-Mod).
\end{Def}

Like the symmetric 2-group case, we have
\begin{Prop}
Let $(F_{\cdot},\lambda_{\cdot}),
(G_{\cdot},\mu_{\cdot}):(\cA_{\cdot},L_{\cdot},\alpha_{\cdot})\rightarrow
(\cB_{\cdot},M_{\cdot},\beta_{\cdot})$ be two morphisms of complexes
of $\cR$-2-modules. If they are 2-chain homotopy, there is an
equivalence $\cH_{n}(F_{\cdot})\backsimeq \cH_{n}(G_{\cdot})$
between induced morphisms.
\end{Prop}

\begin{Lem}
Let $\cR,\cS$ be 2-rings, $(F_{\cdot},\lambda_{\cdot}),
(G_{\cdot},\mu_{\cdot}):(\cA_{\cdot},L_{\cdot},\alpha_{\cdot})\rightarrow
(\cB_{\cdot},M_{\cdot},\beta_{\cdot})$ be two 2-chain homotopy
morphisms of complexes of $\cR$-2-modules and
$T:(\cR$-2-Mod)$\rightarrow(\cS$-2-Mod) be a 2-functor. Then
$T(F_{\cdot},\lambda_{\cdot})$ is 2-chain homotopic to
$T(G_{\cdot},\mu_{\cdot})$ in ($\cS$-2-Mod).
\end{Lem}
%\begin{proof}
%Assume the complex $\cA_{\cdot}$ as the following diagram
%\begin{center}
%\scalebox{0.9}[0.85]{\includegraphics{p13.3.eps}}
%\end{center}
%For $T$:($\cR$-2-Mod)$\rightarrow$($\cS$-2-Mod) is a 2-functor, we
%get a sequence
%\begin{center}
%\scalebox{0.9}[0.85]{\includegraphics{p13.4.eps}}
%\end{center}
%By the definition of 2-functor, the above sequence is a complex in
%($\cS$-2-Mod).
%\end{proof}

\section{Projective Resolution and Derived 2-Functor in ($\cR$-2-Mod)}

In this section we will construct a projective resolution of any
$\cR$-2-module, define the left derived 2-functor and then give the
basic property of this derived 2-functor.
% All the
%results are similar as the symmetric 2-groups case, so we don't
%prove them.

\begin{Def} Let $\cM$ be an $\cR$-2-module.
A projective resolution of $\cM$ in ($\cR$-2-Mod) is a complex of
$\cR$-2-mosules which is relative 2-exact in each point as in the
following diagram
\begin{center}
\scalebox{0.9}[0.85]{\includegraphics{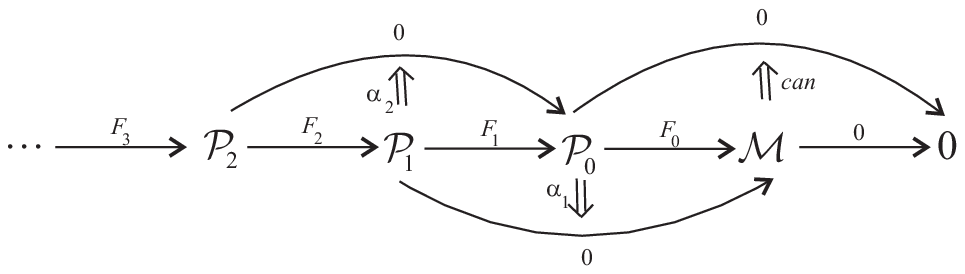}}
\end{center}
with $\cP_{n}(n\geq 0)$ projective objects in ($\cR$-2-Mod). i.e.
the above complex is relative 2-exact in each $\cP_{i}$ and $\cM$.
\end{Def}

\begin{Prop}
Every $\cR$-2-module $\cM$ has a projective resolution in
($\cR$-2-Mod).
\end{Prop}
Sketch of proof. The construction of projective resolution of $\cM$
is similar to symmetric 2-group case.

For $\cM$, there is an essentially surjective morphism
$F_0:\cP_{0}\rightarrow\cM$, with $\cP_{0}$ projective object in
($\cR$-2-Mod)(\cite{14,21}). Then we get a sequence as follows
\begin{center}
\scalebox{0.9}[0.85]{\includegraphics{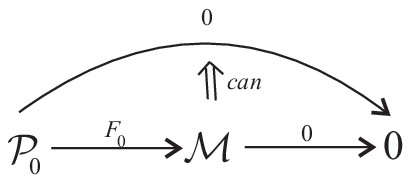}} {\footnotesize S.1.}
\end{center}
where $0:\cM\rightarrow 0$ is the zero morphism\cite{4}in
($\cR$-2-Mod), $0$ is the $\cR$-2-module with only one object and
one morphism., $can$ is the canonical 2-morphism in ($\cR$-2-Mod),
which is given by the identity morphism of only one object of $0$.

From the existence of the relative kernel in ($\cR$-2-Mod), we have
the relative kernel $(Ker(F_0,can),\ e_{(F_0,can)},\
\varepsilon_{(F_0,can)})$
 of the sequence S.1, which is in fact the general kernel
 $(KerF_0,e_{F_0},\varepsilon_{F_0})$ \cite{4}. For $\cR$-2-module $KerF_0$, there exists an essentially surjective morphism
 $G_1:\cP_{1}\rightarrow KerF_0$, with $\cP_{1}$ projective object in
($\cR$-2-Mod)(\cite{14,21}). Let $F_{1}=e_{F_{0}}\circ
 G_{1}:\cP_{1}\rightarrow\cP_{0}$. Then we get the following
 sequence
\begin{center}
\scalebox{0.9}[0.85]{\includegraphics{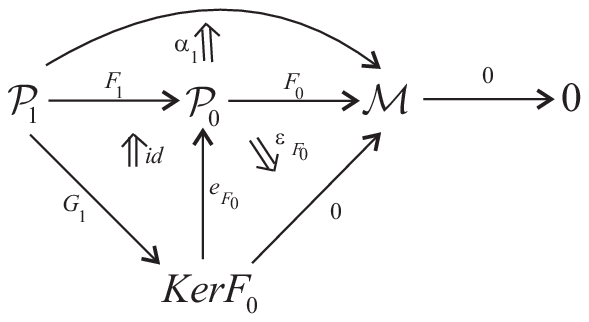}}
%{\footnotesize Fig.1.}
\end{center}
where $\alpha_{1}$ is the composition $F_{0}\circ F_{1}=F_{0}\circ
e_{F_0}\circ G_{1}\Rightarrow 0\circ G_1\Rightarrow 0$ and
compatible with $can$.

Consider the above sequence, there exists the relative kernel
$(Ker(F_1,\alpha_1),e_{(F_1,\alpha_1)},\\\varepsilon_{(F_1,\alpha_1)})$
in (2-SGp). For the $\cR$-2-module $Ker(F_{1},\alpha_1)$, there is
an essentially surjective morphism $G_{2}:\cP_{2}\rightarrow
Ker(F_1,\alpha_1)$, with $\cP_{2}$ projective object in
($\cR$-2-Mod)(\cite{14,21}). Let $F_2=e_{(F_1,\alpha_1)}\circ
G_2:\cP_2\rightarrow\cP_1$. Then we get a sequence
\begin{center}
\scalebox{0.9}[0.85]{\includegraphics{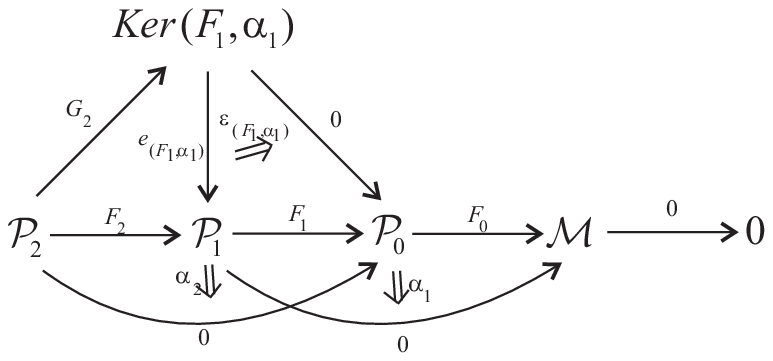}}
%{\footnotesize Fig.1.}
\end{center}
where $\alpha_2$ is the composition $F_{1}\circ F_{2}=F_{1}\circ
e_{(F_1,\alpha_1)}\circ G_{2}\Rightarrow 0\circ G_2\Rightarrow 0$
and compatible with $\alpha_1$.

Using the same method, we get a complex of $\cR$-2-modules
\begin{center}
\scalebox{0.9}[0.85]{\includegraphics{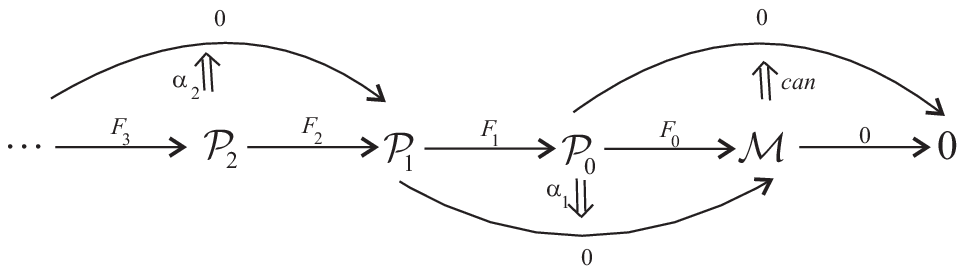}}
%{\footnotesize
%Fig.2.}
\end{center}
The proof of  relative 2-exactness of the sequence is the same as
symmetric 2-group case.

\begin{Thm}
Let $(F_{\cdot}:\cP_{\cdot}\rightarrow \cM,\alpha_{\cdot})$ be a
projective resolution of $\cR$-2-module $\cM$, and $H:\cM\rightarrow
\cN$ a morphism in ($\cR$-2-Mod). Then for any projective resolution
$(G_{\cdot}:\cQ_{\cdot}\rightarrow \cN,\beta_{\cdot})$, there is a
morphism $H_{\cdot}:\cP_{\cdot}\rightarrow \cQ_{\cdot}$ of complexes
in ($\cR$-2-Mod) together with the family of 2-morphisms
$\{\varepsilon_{n}:G_{n}\circ H_{n}\Rightarrow H_{n-1}\circ
F_{n}\}_{n\geq 0}$ as in the following diagram
\begin{center}
\scalebox{0.9}[0.85]{\includegraphics{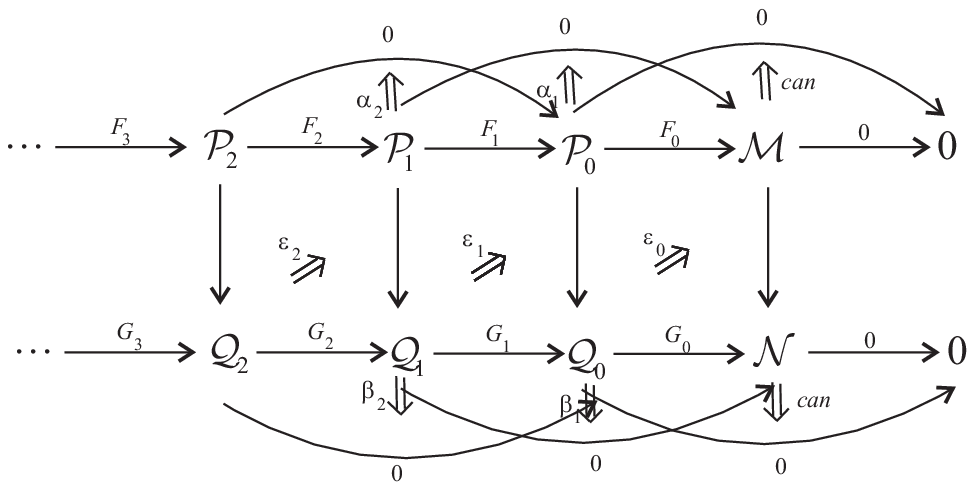}}
\end{center}
If there is another morphism between projective resolutions, they
are 2-chain homotopy.
\end{Thm}
 The proof of this Theorem is also similar to symmetric
2-group case in \cite{20}. The difference is that the existence of
1-morphisms and 2-morphisms is from the properties of projective
$\cR$-2-modules in ($\cR$-2-Mod).

\begin{Def}Let $\cR,\cS$ be 2-rings.
An additive 2-functor(\cite{2}) $T$:
($\cR$-2-Mod)$\rightarrow$($\cS$-2-Mod) is called right relative
2-exact if the relative 2-exactness of
\begin{center}
\scalebox{0.9}[0.85]{\includegraphics{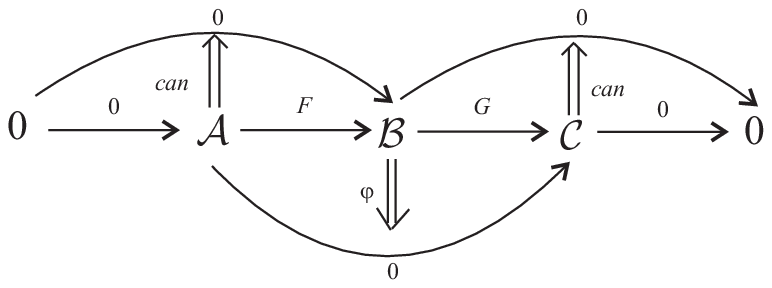}}
\end{center}
in $\cA,\cB$ and $\cC$ implies relative 2-exactness of
\begin{center}
\scalebox{0.9}[0.85]{\includegraphics{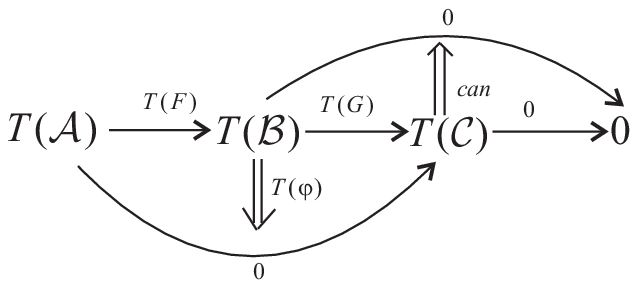}}
\end{center}
in $T(\cB)$ and $T(\cC)$.
\end{Def}
The left relative 2-exact 2-functor can be defined dually.

By Remark 2 and Proposition 1, Theorem 1, there is
\begin{cor}
Let $\cR,\cS$ be 2-rings, and $T:
$($\cR$-2-Mod)$\rightarrow$($\cS$-2-Mod) be an additive 2-functor,
and $\cA$ be any object of ($\cR$-2-Mod). For two projective
resolutions $\cP_{\cdot},\ \cQ_{\cdot}$ of $\cA$, there is an
equivalence between homology $\cS$-2-modules
$\cH_{\cdot}(T(\cP_{\cdot}))$ and $\cH_{\cdot}(T(\cQ_{\cdot}))$.
\end{cor}

Let $T$: ($\cR$-2-Mod)$\rightarrow$($\cS$-2-Mod) be an additive 2-functor.
There is a 2-functor
\begin{align*}
&\cL_{i}T:(\cR\textrm{-2-Mod)}\rightarrow (\cS\textrm{-2-Mod)}\\
&\hspace{2.7cm}\cA\mapsto \cL_{i}T(\cA),\\
&\hspace{1.7cm}\cA\xrightarrow[]{F}\cB\mapsto
\cL_{i}T(\cA)\xrightarrow[]{\cL_{i}T(F)} \cL_{i}T\cB),
\end{align*}
where $\cL_{i}T(\cA)$ is defined by $\cH_{i}(T(\cP_{\cdot}))$, and
$\cP_{\cdot}$ is the projective resolution of $\cA$. $\cL_{i}T$ is a
well-defined 2-functor from the properties of additive 2-functor and
Corollary 1.
\begin{cor}
Let $T$: ($\cR$-2-Mod)$\rightarrow$($\cS$-2-Mod) be a right relative
2-exact 2-functor, and $\cA$ be a projective object in
($\cR$-2-Mod). Then $\cL_{i}T(\cA)=0$ for $i\neq0$.
\end{cor}
%we can construct the left derived 2-functor $\cL_{i}T(i\geq 0)$ of
%$T$ as follows. If $\cA$ is an object of (2-SGp), choose a
%projective resolution $\cP_{\cdot}\rightarrow\cA$ and define
%$$
%\cL_{i}T(\cA)=\cH_{i}(T(\cP_{\cdot})).
%$$
%If $H:\cA\rightarrow\cB $ is any 1-morphism in (2-SGp), by Theorem
%1, there is a morphism of projective resolutions of $\cA$ and $\cB$,
%and then there is a morphism of $\cL_{i}T(\cA)$ and $\cL_{i}T(\cB)$.

%\begin{cor} If $\cA$ is a projective object in (2-SGp), then
%$\cL_{i}(T(\cA))=0$ for $i\neq 0$.
%\end{cor}
The following is a basic property of derived functors.
\begin{Thm}Let $T$: ($\cR$-2-Mod)$\rightarrow$($\cS$-2-Mod) be a right relative
2-exact 2-functor. Then the left derived 2-functor $\cL_{*}T$ takes
%form a homological functor, i.e.
the sequence of $\cR$-2-modules
\begin{center}
\scalebox{0.9}[0.85]{\includegraphics{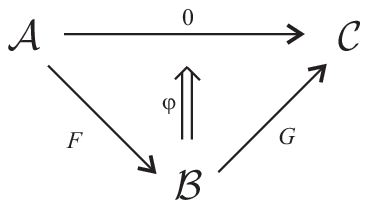}}
\end{center}
which is relative 2-exact in $\cA,\ \cB,\ \cC$ to a long sequence
2-exact(similar \cite{1,6})in each point
\begin{center}
\scalebox{0.9}[0.85]{\includegraphics{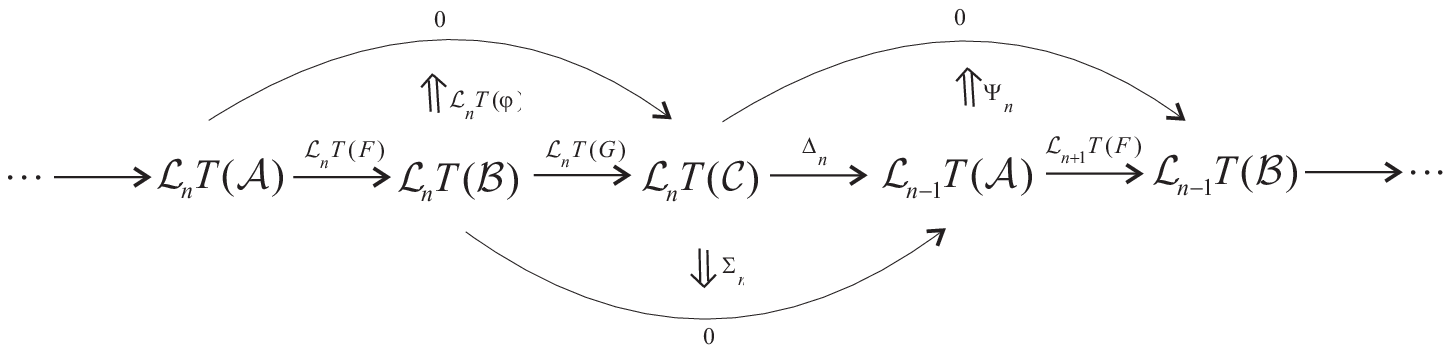}}
\end{center}
\end{Thm}

In order to prove this theorem, we need the following Lemmas.

Similar to the proofs of symmetric 2-group case. We have
\begin{Lem}
Let $\cP$ and $\cQ$ be projective objects in ($\cR$-2-Mod). Then the
product category $\cP\times\cQ$ is a projective object in
($\cR$-2-Mod).
\end{Lem}

\begin{Lem}
Let $(F,\varphi,G):\cA\rightarrow\cB\rightarrow\cC$ be an extension
of $\cR$-2-modules in ($\cR$-2-Mod)(similar to \cite{1,11}),
$(\cP_{\cdot},L_{\cdot},\alpha_{\cdot})$
$(\cQ_{\cdot},N_{\cdot},\beta_{\cdot})$ be projective resolutions of
$\cA$ and $\cC$, respectively. Then there is a projective resolution
$(\cK_{\cdot},M_{\cdot},\varphi_{\cdot})$ of $\cB$, such that
$\cP_{\cdot}\rightarrow \cK_{\cdot}\rightarrow\cQ_{\cdot}$ forms an
extension of complexes in ($\cR$-2-Mod).
\end{Lem}

By the universal property of (bi)product of $\cR$-2-modules and the
property of additive 2-functor(\cite{2,4}). We get
\begin{Lem}
Let $T$: ($\cR$-2-Mod)$\rightarrow$($\cS$-2-Mod) be an additive
2-functor, and $\cA,\ \cB$ be objects in ($\cR$-2-Mod). Then there
is an equivalence between $T(\cA\times\cB)$ and $T(\cA)\times
T(\cB)$ in ($\cS$-2-Mod).
\end{Lem}

Proof of Theorem 2. For $\cR$-2-modules $\cA$ and $\cC$, choose
projective resolutions $\cP_{\cdot}\rightarrow\cA$ and
$\cQ_{\cdot}\rightarrow\cC$. By Lemma 2 and Lemma 3, there is a
projective resolution $\cP_{\cdot}\times\cQ_{\cdot}\rightarrow \cB$
fitting into an extension
$\cP_{\cdot}\xrightarrow[]{i_{\cdot}}\cP_{\cdot}\times\cQ_{\cdot}\xrightarrow[]{p_{\cdot}}\cQ_{\cdot}$
of projective complexes in (2-SGp)(\cite{1}). By Lemma 4, we obtain
a complexes of extension
$$T(\cP_{\cdot})\xrightarrow[]{T(i_{\cdot})}T(\cP_{\cdot}\times\cQ_{\cdot})\xrightarrow[]{T(p_{\cdot})}T(\cQ_{\cdot}).$$

Similar to Theorem 4.2 in \cite{11}, the long sequence
\begin{center}
\scalebox{0.9}[0.85]{\includegraphics{p28.eps}}
\end{center}
is 2-exact in each point.

\section*{Acknowledgements.} We thank Prof. Zhang-Ju LIU, Prof. Yun-He SHENG for providing
us the ideas of higher dimensional category theory. We also thank
Prof. Ke WU and Prof. Shi-Kun WANG for useful discussions.

\noindent Fang Huang, Shao-Han Chen, Wei Chen\\
Department of Mathematics\\
 South China University of
 Technology\\
 Guangzhou 510641, P. R. China

\noindent Zhu-Jun Zheng\\
Department of Mathematics\\
 South China University of
 Technology\\
 Guangzhou 510641, P. R. China \\
 and\\
Institute of Mathematics\\
Henan University\\  Kaifeng 475001, P. R.
China\\
E-mail: zhengzj@scut.edu.cn

\end{document}